\documentclass[graybox]{svmult}

\usepackage{type1cm}        % activate if the above 3 fonts are
\usepackage{makeidx}         % allows index generation
\usepackage{graphicx}        % standard LaTeX graphics tool
\usepackage{multicol}        % used for the two-column index
\usepackage[bottom]{footmisc}% places footnotes at page bottom

\usepackage{newtxtext}       % 
\usepackage{newtxmath}       % selects Times Roman as basic font
\usepackage{amsmath}

\makeindex             % used for the subject index
\begin{document}

% \title*{An IGA Framework for PDE-Based Planar Parameterization with Arbitrary Interface Continuity}
\title*{An IGA Framework for PDE-Based Planar Parameterization on Convex Multipatch Domains}
%\title*{PDE-based Planar Parameterization with Support for Locally Reduced Smoothness and Convex Multipatch Domains}
% \titlerunning{Planar Parameterization with Arbitrary Interface Continuity}
\titlerunning{Planar Parameterization on Connvex Multipatch Domains}
\author{Jochen Hinz, Matthias M\"oller, Cornelis Vuik}
%\institute{Jochen Hinz\at Delft Institute of Applied Mathematics, Mourik Broekmanweg 6, 2628XE Delft, the Netherlands, \email{j.p.hinz@tudelft.nl}
%\and Matthias M\"oller \at Delft Institute of Applied Mathematics, Mourik Broekmanweg 6, 2628XE Delft, the Netherlands, \email{m.moller@tudelft.nl}
%\and Cornelis Vuik \at Delft Institute of Applied Mathematics, Mourik Broekmanweg 6, 2628XE Delft, the Netherlands, \email{c.vuik@tudelft.nl}}
\institute{Jochen Hinz, \email{j.p.hinz@tudelft.nl}
\and Matthias M\"oller, \email{m.moller@tudelft.nl}
\and Cornelis Vuik, \email{c.vuik@tudelft.nl} \at Delft Institute of Applied Mathematics, Mourik Broekmanweg 6, 2628XE Delft, the Netherlands}
\maketitle

\abstract{The first step towards applying isogeometric analysis techniques to solve PDE problems on a given domain consists in generating an analysis-suitable mapping operator between parametric and physical domains with one or several patches from no more than a description of the boundary contours of the physical domain. A subclass of the multitude of the available parameterization algorithms are those based on the principles of \emph{Elliptic Grid Generation} (EGG) which, in their most basic form, attempt to approximate a mapping operator whose inverse is composed of harmonic functions. The main challenge lies in finding a formulation of the problem that is suitable for a computational approach and a common strategy is to approximate the mapping operator by means of solving a PDE-problem. PDE-based EGG is well-established in classical meshing and first generalization attempts to spline-based descriptions (as is mandatory in IgA) have been made. Unfortunately, all of the practically viable PDE-based approaches impose certain requirements on the employed spline-basis, in particular global $C^{\geq 1}$-continuity. \\
This paper discusses an EGG-algorithm for the generation of planar parameterizations with locally reduced smoothness (i.e., with support for locally only $C^0$-continuous bases). A major use case of the proposed algorithm is that of multipatch parameterizations, made possible by the support of $C^0$-continuities. This paper proposes a specially-taylored solution algorithm that exploits many characteristics of the PDE-problem and is suitable for large-scale applications. It is discussed for the single-patch case before generalizing its concepts to multipatch settings. This paper is concluded with three numerical experiments and a discussion of the results.}

\section{Introduction}
\label{sect:Introduction}
The automatic generation of analysis-suitable planar parameterizations for IgA-based numerical simulations is a difficult, yet important problem in the field of isogeometric analysis, since generally no more than a description of the boundary contours is available. The main challenge lies in the generation of a folding-free (i.e., bijective) parameterization with numerically favorable properties such as orthogonal isolines and a large degree of parametric \emph{smoothness}. Furthermore, a practical algorithm should be computationally inexpensive, and, if possible, exhibit little sensitivity to small perturbations in the boundary contour description.\\
Let $\Omega$ denote the target geometry and $\hat{\Omega}$ the parametric domain.  Furthermore, let $\mathbf{x}: \hat{\Omega} \rightarrow \Omega$ denote the mapping operator that we attempt to build from the linear span of the B-Spline basis $\Sigma = \{ w_1, w_2, \ldots, w_N \}$, where $\mathbf{x} \vert_{\partial \hat{\Omega} } = \partial \Omega$ is known. Note that $\mathbf{x}$ is of the form:
\begin{align}
\label{eq:form_x}
\mathbf{x}(\xi, \eta) = \sum_{i \in \mathcal{I}_{\text{boundary}}} \mathbf{c}_i w_i(\xi, \eta) + \sum_{j \in \mathcal{I}_{\text{inner}}} \mathbf{c}_j w_j(\xi, \eta),
\end{align}
where $\mathcal{I}_{\text{inner}}$ and $\mathcal{I}_{\text{boundary}}$ denote the index set of the vanishing and nonvanishing basis functions on $\partial \hat{\Omega}$, respectively. Formally, $\mathcal{I}_{\text{boundary}} \cap \mathcal{I}_{\text{inner}} = \emptyset$ and $\mathcal{I}_{\text{boundary}} \cup \mathcal{I}_{\text{inner}} = \{1, \ldots, N\}$. With this, the objective of all parameterization algorithms is to properly select the inner control points $\mathbf{c}_j$, while the boundary control points $\mathbf{c}_i$ are known from the boundary contours and typically held fixed. \\
In \cite{gravesen2012planar}, Gravesen et al. study planar parameterization techniques based on the constrained minimization of a quality functional over the inner control points. To avoid self-intersections, a nonlinear and nonconvex sufficient condition for $\det J > 0$, where $J$ denotes the Jacobian of the mapping, is added as a constraint. The numerical quality of the resulting parameterization depends on the choice of the employed cost functional and the characteristic properties of $\Omega$. While this approach is not guaranteed to yield acceptable results for all types of geometries (see section \ref{sect:Numerical_Experiments}), it is known to yield good results in a wide range of applications with proper parameter tuning. A drawback is the relatively large number of required iterations (typically $\sim 30$) and the need to find an initial guess that satisfies the constraints (for which another optimization problem has to be solved first). The proposed minimization is tackled with a black-box nonlinear optimizer (IPOPT \cite{biegler2009large}) that comes with all the drawbacks of nonlinear optimization such as the danger of getting stuck in local minima. \\
Another class of parameterization methods suitable for nontrivial geometries are PDE-based, most notably, the class of methods based on the principles of \emph{elliptic grid generation} (EGG).  Methods based on EGG attempt to generate a mapping $\mathbf{x}: \hat{\Omega} \rightarrow \Omega$ such that the components of $\mathbf{x}^{-1}: \Omega \rightarrow \hat{\Omega}$ are harmonic functions on $\Omega$. For this, a nonlinear partial differential equation is imposed on $\mathbf{x}$, which takes the form
\begin{align}
\label{eq:EGG_not_discretized}
\mathcal{L}(\mathbf{x}) = g_{22} \mathbf{x}_{\xi \xi} - 2 g_{12} \mathbf{x}_{\xi \eta} + g_{11} \mathbf{x}_{\eta \eta} = 0, \quad \text{s.t. } \mathbf{x} \vert_{\partial \hat{\Omega}} = \partial \Omega,
\end{align}
with
\begin{align}
g_{11}(\mathbf{x}) &= \mathbf{x}_\xi \cdot \mathbf{x}_\xi, \nonumber \\
g_{12}(\mathbf{x}) &= \mathbf{x}_\xi \cdot \mathbf{x}_\eta, \nonumber \\
g_{22}(\mathbf{x}) &= \mathbf{x}_\eta \cdot \mathbf{x}_\eta
\end{align}
 being the entries of the metric tensor of the mapping (which are nonlinear functions of $\mathbf{x}$). Under certain assumptions of the boundary contour regularity and assuming that $\hat{\Omega}$ is convex, it can be shown that the exact solution of (\ref{eq:EGG_not_discretized}) is bijective, justifying a numerical approximation for the purpose of generating a geometry description \cite{azarenok2009generation}. \\
EGG has been an established approach in classical meshing for decades and first attempts to apply it to spline-based geometry descriptions were made in \cite{manke1989tensor}, where the equations are approximately solved by a collocation at the abscissae of a Gaussian quadrature scheme with cubic Hermite-splines. In \cite{lamby2007elliptic}, the collocation takes place at the Greville-abscissae and the resulting nonlinear equations are solved using a Picard-based iterative scheme, allowing for  a wider range of spline-bases. However, as a downside, the consistency order of Greville-based collocation is not optimal. In \cite{hinz2018elliptic}, the equations are discretized with a Galerkin approach and a Newton-based iterative approach is employed for the resulting root-finding problem, allowing for $C^{ \geq 1}$-continuous bases. Numerical convergence is accelerated by generating good initial guesses utilizing multigrid-techniques and convergence is typically achieved within $4$ (unconstrained) nonlinear iterations. \\
Unfortunately, none of the aforementioned approaches allow for spline-bases with locally reduced smoothness, limiting their usefullness in practice, since in certain applications $C^0$-continuities are desirable or unavoidable, notably in multipatch parameterizations or when $\partial \Omega$ is build from a spline-basis with (one or more) $p$-fold internal knot repetitions (where $p$ refers to the polynomial order of the spline-basis used). To allow for $C^0$-continuities, one may instead minimize the \emph{Winslow-functional} \cite{winslow1981adaptive} (whose global minimizer is equal to the exact solution of (\ref{eq:EGG_not_discretized})). Unfortunately, this leads to a formulation in which the Jacobian determinant appears in the denominator, which is why an iterative solution scheme has to be initialized with a bijective initial guess in order to avoid division by zero, restricting it to use cases in which a bijective initial guess is available. \\
Motivated by our striving for a computationally inexpensive parameterization algorithm that does not have to be initialized by a bijective initial guess and allows for spline-bases with arbitrary continuity properties, in this paper, we augment the discretization proposed in \cite{hinz2018elliptic} with auxilliary variables, leading to a mixed-FEM type problem. To allow for its application to large-scale problems, we present a solution strategy that tackles the resulting nonlinear root-finding problem with a Newton-Krylov-based \cite{knoll2004jacobian} Jacobian-free iterative approach that only operates on the nonlinear part (corresponding to the primary, not auxilliary variables) of the equation. Besides single-patch problems, we will address potential use cases of the algorithm in multipatch settings (in particular with extraordinary vertices), made possible by the support of $C^0$-continuous spline bases. We conclude this paper with a number of example-parameterizations and a discussion of the results.

\section{Problem Formulation}
In \cite{hinz2018elliptic}, the following discretization of the governing equations (see equation (\ref{eq:EGG_not_discretized})) is proposed:
\begin{align}
\label{eq:main_C1_unscaled}
& \text{find } \mathbf{x} \in [\operatorname{span} \Sigma]^2 \text{ s.t.} \nonumber \\
& \left \{ \begin{array}{l} \forall \boldsymbol{\sigma}_i \in [\Sigma_0]^2: \quad \int_{\hat{\Omega}} \boldsymbol{\sigma}_i \cdot \mathcal{L}(\mathbf{x}) \mathrm{d} \boldsymbol{\xi} = 0 \\
			 \mathbf{x}\vert_{\partial \hat{\Omega}} = \partial \Omega \end{array} \right.,
\end{align}
where $\Sigma_0 \equiv \{ w_i \in \Sigma \enskip \vert \enskip w_i \vert_{\partial \hat{\Omega}} = 0 \}$. \\
Similarly, \cite{hinz2018spline} introduces a scaled version of (\ref{eq:main_C1_unscaled}), namely:
\begin{align}
\label{eq:main_C1_scaled}
& \text{find } \mathbf{x} \in [\operatorname{span} \Sigma]^2 \text{ s.t.} \nonumber \\
& \left \{ \begin{array}{l} \forall \boldsymbol{\sigma}_i \in [\Sigma_0]^2: \quad \int_{\hat{\Omega}} \boldsymbol{\sigma}_i \cdot \tilde{\mathcal{L}}(\mathbf{x}) \mathrm{d} \boldsymbol{\xi} = 0 \\ 
\mathbf{x}\vert_{\partial \hat{\Omega}} = \partial \Omega \end{array} \right.,
\end{align}
where
\begin{align}
\label{eq:L_scaled}
\tilde{\mathcal{L}}( \mathbf{x} ) = \frac{\mathcal{L}(\mathbf{x})}{ \underbrace{g_{11} + g_{22}}_{ \geq 0 } + \underbrace{\mu}_{> 0} }.
\end{align}
Here, $\mu > 0$ is a small positive parameter that is usually taken to be $\mu = 10^{-4}$. \\
The motivation to solve (\ref{eq:main_C1_scaled}) rather than (\ref{eq:main_C1_unscaled}) is based on the observation that numerical root-finding algorithms typically converge faster in this case and that a suitable convergence criterion is less geometry-dependent. Note that the scaling is allowed because the exact solution is unchanged. Therefore, we base our reformulation of the problem on (\ref{eq:main_C1_scaled}). \\
In order to reduce the highest-order derivatives from two to one, we introduce a new operator in which we replace second order derivatives in $\mathbf{x}$ by the first order derivatives of $\mathbf{u}$ and $\mathbf{v}$, respectively:
\begin{align}
\label{eq:U_xi_eta}
\mathcal{U}( \mathbf{u}, \mathbf{v}, \mathbf{x} ) & = \frac{ g_{22} \mathbf{u}_{\xi} - g_{12} \mathbf{u}_{\eta} - g_{12} \mathbf{v}_{\xi} + g_{11} \mathbf{v}_{\eta} }{ g_{11} + g_{22} + \mu }.
\end{align}
Where $\mathcal{U}$ satisfies
\begin{align}
\tilde{\mathcal{L}}( \mathbf{x} ) = \mathcal{U}( \mathbf{x}_\xi, \mathbf{x}_\eta, \mathbf{x} ).
\end{align}
A possible reformulation of (\ref{eq:main_C1_scaled}) with auxilliary variables now reads:
\begin{align}
\label{eq:main_xi_eta}
& \text{find } (\mathbf{u}, \mathbf{v}, \mathbf{x})^T \in [ \operatorname{span} \bar{\Sigma} ]^4 \times [ \operatorname{span} \Sigma ]^2 \text{ s.t.} \nonumber \\
& \left \{ \begin{array}{l} \forall \boldsymbol{\sigma}_i \in [\bar{\Sigma}]^4 \times [\Sigma_0]^2: \quad \int_{\hat{\Omega}} \boldsymbol{\sigma}_i \cdot \begin{pmatrix} \mathbf{u} - \mathbf{x}_\xi \\ \mathbf{v} - \mathbf{x}_\eta \\ \mathcal{U}(\mathbf{u}, \mathbf{v}, \mathbf{x}) \end{pmatrix} \mathrm{d} \boldsymbol{\xi} = 0 \\ \mathbf{x}\vert_{\partial \hat{\Omega}} = \partial \Omega \end{array} \right.,
\end{align}
where $\bar{\Sigma} = \{ \bar{w}_1, \ldots, \bar{w}_{\bar{N}} \}$ denotes the basis that is used for the auxilliary variables. \\
Note that the choice of (\ref{eq:U_xi_eta}) is not unique. Here, we have chosen to divide $\mathbf{x}_{\xi \eta}$ equally among $\mathbf{u}_\eta$ and $\mathbf{v}_\xi$. In general, any combination
\begin{align}
\mathbf{x}_{\xi \eta} \rightarrow \chi \mathbf{u}_\eta + (1 - \chi) \mathbf{v}_\xi,
\end{align}
is valid. Note that since the $g_{ij}$ are functions of $\mathbf{x}_\xi$ and $\mathbf{x}_\eta$, further possible variants are acquired by substituting $\mathbf{u}, \mathbf{v}$ in the $g_{ij}$. \\
System (\ref{eq:main_xi_eta}) now constitutes a discretization of (\ref{eq:EGG_not_discretized}) that allows for only $C^0$-continuous bases at the expense of increasing the problem size from $2 \vert \mathcal{I}_{\text{inner}} \vert$ to $2 \vert \mathcal{I}_{\text{inner}} \vert + 4 \vert \bar{\Sigma} \vert$, where, as before, $\mathcal{I}_{\text{inner}}$ refers to the index set of inner control points. \\
Let us remark that in certain settings, it suffices to invoke auxilliary variables in one coordinate-direction only. A possible problem formulation for the $\xi$-direction reads:
\begin{align}
\label{eq:main_xi}
& \text{find } (\mathbf{u}, \mathbf{x})^T \in [\operatorname{span} \bar{\Sigma}]^2 \times [\operatorname{span} \Sigma]^2 \text{ s.t.} \nonumber \\
& \left \{ \begin{array}{l} \forall \boldsymbol{\sigma}_i \in [\bar{\Sigma}]^2 \times [\Sigma_0]^2: \quad \int_{\hat{\Omega}} \boldsymbol{\sigma}_i \cdot \begin{pmatrix} \mathbf{u} - \mathbf{x}_\xi \\ \mathcal{U}^\xi(\mathbf{u}, \mathbf{x}) \end{pmatrix} \mathrm{d} \boldsymbol{\xi} = 0 \\ \mathbf{x}\vert_{\partial \hat{\Omega}} = \partial \Omega \end{array} \right.,
\end{align}
with (for instance)
\begin{align}
\label{eq:U_xi}
\mathcal{U}^\xi (\mathbf{u}, \mathbf{x} ) & = \frac{ g_{22} \mathbf{u}_{\xi} - g_{12} \mathbf{u}_{\eta} - g_{12} \mathbf{x}_{\xi \eta} + g_{11} \mathbf{x}_{\eta \eta} }{ g_{11} + g_{22} + \mu }.
\end{align}
And similarly for the $\eta$-direction. \\
The above approach is useful if $C^0$-continuities are only required in a single coordinate-direction so that the total number of degrees of freedom (DOFs) can be reduced.

\section{Solution Strategy}
Systems (\ref{eq:main_xi_eta}) and (\ref{eq:main_xi}) are nonlinear and have to be solved with an iterative algorithm. We will discuss a solution algorithm that is losely based on the Newton-approach proposed in \cite{hinz2018elliptic}. However, we tweak it in order to reduce computational costs and memory requirements by exploiting many characteristics of the problem at hand. First, we discuss the case in which $\hat{\Omega}$ is given by a single patch, after which we generalize our solution strategy to multipatch-settings (in particular with topologies that contain extraordinary vertices).

\subsection{Single Patch Paramtererizations}
\label{subsect:Solution_Single}
With $\mathbf{x} = \mathbf{x}[ \mathbf{c} ]$, where $\mathbf{c}$ is a vector containing the $\mathbf{c}_j$ in (\ref{eq:form_x}) (while freezing the $\mathbf{c}_i$ that follow from the boundary condition) and $(\mathbf{u}, \mathbf{v})^T = (\mathbf{u}, \mathbf{v})^T[ \mathbf{d} ]$, where $\mathbf{d} = (\mathbf{d}^u, \mathbf{d}^v)^T$ is a vector containing $\mathbf{d}_i^{u}$ and $\mathbf{d}_i^{v}$ in
\begin{align}
\mathbf{u}[\mathbf{d}^u] & = \sum_i \mathbf{d}_i^u \bar{w_i}, \nonumber \\
\mathbf{v}[\mathbf{d}^v] & = \sum_i \mathbf{d}_i^v \bar{w_i},
\end{align}
we can reinterpret (\ref{eq:main_xi_eta}) as a problem in $\mathbf{c}$ and $\mathbf{d}$. It has a residual vector of the form
\begin{align}
\label{eq:R_ux}
\mathbf{R}( \mathbf{d}, \mathbf{c} ) = \begin{pmatrix} R_L ( \mathbf{d}, \mathbf{c} ) \\ R_N( \mathbf{d}, \mathbf{c} )  \end{pmatrix},
\end{align}
where $R_L$ refers to the linear part in (\ref{eq:main_xi_eta}) (the projection of the auxilliary variables onto $\mathbf{x}_\xi$ and $\mathbf{x}_\eta$) and $R_N$ to the nonlinear (the part involving the operator $\mathcal{U}(\mathbf{u}, \mathbf{v}, \mathbf{x})$). \\
The Newton-approach from \cite{hinz2018elliptic} requires the assembly of the Jacobian
\begin{align}
J_R = 
\begingroup
\setlength\arraycolsep{2pt}
\renewcommand{\arraystretch}{1.5}
\begin{pmatrix} \frac{ \partial R_L }{ \partial \mathbf{d} } & \frac{ \partial R_L }{ \partial \mathbf{c} } \\
				    \frac{ \partial R_N }{ \partial \mathbf{d} } & \frac{ \partial R_N }{ \partial \mathbf{c} } \end{pmatrix}
\endgroup
\equiv \begin{pmatrix} A & B \\				   																										C & D \end{pmatrix}
\end{align} 
 of (\ref{eq:main_xi_eta}) at every Newton-iteration. The matrices $A$ and $B$, corresponding to the linear part in (\ref{eq:main_xi_eta}), are not a function of $\mathbf{c}$ and $\mathbf{d}$ and thus have to be assembled only once. In fact, $A$ is block-diagonal with blocks given by the parametric mass matrix $\bar{M}$ over the auxilliary basis $\bar{\Sigma} = \{ \bar{w}_1, \ldots, \bar{w}_{\bar{N}} \}$ with entries
\begin{align}
\bar{M}_{ij} = \int_{\hat{\Omega}} \bar{w}_i \bar{w}_j \mathrm{d} \boldsymbol{\xi},
\end{align}
while $B$ is block-diagonal with blocks whose columns are given by a subset of the colums of the matrices $\bar{M}^\xi$ and $\bar{M}^\eta$ with entries
\begin{align}
\bar{M}^\xi_{ij} = \int_{\hat{\Omega}} \bar{w}_i w_{j\xi} \mathrm{d} \boldsymbol{\xi}
\end{align}
and
\begin{align}
\bar{M}^\eta_{ij} = \int_{\hat{\Omega}} \bar{w}_i w_{j\eta} \mathrm{d} \boldsymbol{\xi}.
\end{align}
%Hence, both $A$ and $B$ are separable matrices and can therefore be assembled inexpensively from a kronecker-product of smaller matrices containing univariate integrals only. To further improve the efficiency, one may choose to never carry out the kronecker product explicitly and only store the univariate constituents of $A$ and $B$. \\
For given $\mathbf{c}$ and $\mathbf{d}$, the Newton search-direction is computed from a system of the form
\begin{align}
\label{eq:Jacobian_full}
\begin{pmatrix} A & B \\ C & D \end{pmatrix} \begin{pmatrix} \Delta \mathbf{d} \\ \Delta \mathbf{c} \end{pmatrix} & = \begin{pmatrix} \mathbf{a} \\ \mathbf{b} \end{pmatrix},
\end{align}
where $C = C(\mathbf{d}, \mathbf{c})$ and $D = D(\mathbf{d}, \mathbf{c})$ are, unlike $A$ and $B$, not constant and have to be reassembled in each iteration. We form the Schur-complement of $A$, in order to yield an equation for $\Delta \mathbf{c}$ only, namely:
\begin{align}
\label{eq:Schur_equation}
(\underbrace{D - C A^{-1} B}_{\tilde{D}}) \Delta \mathbf{c} = \mathbf{b} - C A^{-1} \mathbf{a}.
\end{align}
In order to avoid the computationally expensive assembly of $C$ and $D$, we solve (\ref{eq:Schur_equation}) with a Newton-Krylov \cite{knoll2004jacobian} algorithm which only requires the evaluation of vector products $\tilde{D} \mathbf{s}$, which can be approximated with finite differences rather than explicit assembly of $C$ and $D$. Since
\begin{align}
C \mathbf{s}_1 + D \mathbf{s}_2 = \frac{ R_N \left( \mathbf{d} + \epsilon \mathbf{s}_1, \mathbf{c} + \epsilon \mathbf{s}_2 \right) - R_N(\mathbf{d}, \mathbf{c}) }{ \epsilon } + \mathcal{O}(\epsilon),
\end{align}
we have
\begin{align}
\label{eq:D_vec_approx}
\tilde{D} \mathbf{s} & \simeq \frac{ R_N( \mathbf{d} - \epsilon A^{-1} B \mathbf{s}, \mathbf{c} + \epsilon \mathbf{s} )  - R_N( \mathbf{d}, \mathbf{c} ) }{ \epsilon },
\end{align}
and
\begin{align}
\label{eq:C_vec_approx}
C A^{-1} \mathbf{a} \simeq \frac{ R_N( \mathbf{d} + \epsilon A^{-1} \mathbf{a}, \mathbf{c} ) - R_N( \mathbf{d}, \mathbf{c} ) }{ \epsilon },
\end{align}
for $\epsilon$ small. The optimal choice of $\epsilon$ is discussed in \cite{knoll2004jacobian}.\\
We compute products of the form $\mathbf{q} = A^{-1} \mathbf{t}$ from the solution of the system $A \mathbf{q} = \mathbf{t}$, which has for $\mathbf{t} = B \mathbf{s}$ (see equation (\ref{eq:D_vec_approx})) and $\mathbf{t} = \mathbf{a}$ (see equation (\ref{eq:C_vec_approx})) the form of a (separable) $L_2$-projection. Let
\begin{align}
\mathbf{x}^0 [\mathbf{c}] = \sum_{j \in \mathcal{I}_\text{inner} } \mathbf{c}_j w_j.
\end{align}
Product $\mathbf{q} = A^{-1} B \mathbf{s}$ satisfies
\begin{align}
\mathbf{q} = ( \mathbf{q}^u, \mathbf{q}^v )^T = \underset{ (\tilde{\mathbf{q}}^u, \tilde{\mathbf{q}}^v) }{ \operatorname{argmin} } \frac{1}{2}\int_{\hat{\Omega}} \left \| \begin{bmatrix} \mathbf{u}[\tilde{\mathbf{q}}^u] \\ \mathbf{v}[\tilde{\mathbf{q}}^v] \end{bmatrix}  - \begin{bmatrix} \mathbf{x}_\xi^0[ \mathbf{s} ] \\ \mathbf{x}_\eta^0[ \mathbf{s} ]  \end{bmatrix}  \right \|^2 \mathrm{d} \boldsymbol{\xi},
\end{align}
and similarly for $\mathbf{q} = A^{-1} \mathbf{a}$. \\
 As such, $A$ is block-diagonal and composed of separable mass matrices $\bar{M} = \bar{m}_\xi \otimes \bar{m}_\eta$
\begin{align}
A = \begin{pmatrix} \bar{m}_\xi \otimes \bar{m}_\eta &  & \\ & \ddots & \\ & & \bar{m}_\xi \otimes \bar{m}_\eta \end{pmatrix},
\end{align}
where $\bar{m}_\xi$ and $\bar{m}_\eta$ refer to the univariate mass matrices resulting from the tensor-product structure of $\bar{\Sigma}$. Therefore, we have
\begin{align}
A^{-1} = \begin{pmatrix} ( \bar{m}_\xi^{-1} ) \otimes ( \bar{m}_\eta^{-1} ) &  & \\ & \ddots & \\ & & ( \bar{m}_\xi^{-1} ) \otimes ( \bar{m}_\eta^{-1} ) \end{pmatrix}.
\end{align}
We follow the methodology from \cite{gao2014fast}, where a computationally inexpensive inversion of this $2D$ mass matrix is achieved by repeated inversion with the $1D$ mass matrices $\bar{m}_{\xi}$ and $\bar{m}_{\eta}$. Here, we do direct inversion of the $1D$ mass matrices by computing their Cholesky-decompositions \cite{seiler1989numerical}. An inversion can be done in only $\mathcal{O}(\bar{N})$ arithmetic operations and Cholesky-decompositions have to be formed only once, thanks to the fact that $A$ is constant. \\
After solving (\ref{eq:Schur_equation}), $\Delta \mathbf{d}$ is found by solving
\begin{align}
A \Delta \mathbf{d} = \mathbf{a} - B \Delta \mathbf{c}.
\end{align}
Upon completion, the vector $\mathbf{n} \equiv ( \Delta \mathbf{d}, \Delta \mathbf{c} )^T$ constitutes the Newton search-direction. We update the current iterate $( \mathbf{d}, \mathbf{c} )^T$ by adding $\nu \mathbf{n}$, where the optimal value of $\nu \in (0, 1]$ is estimated through a line-search routine. Above steps are repeated until the norm of $\mathbf{n}$ is negligibly small. Upon completion, we extract the $\mathbf{c}$-component from the resulting solution vector which contains the inner control points of the mapping operator $\mathbf{x}$, while the $\mathbf{d}$-component serves no further purpose and can be discarded. \\
It should be noted that a single matrix-vector product $\tilde{D} \mathbf{s}$ is slightly more expensive than, for instance, $D \mathbf{s}$, due to the requirement to invert $A$. However, thanks to the separable nature of $A$, the costs in (\ref{eq:D_vec_approx}) are dominated by function evaluations in $R_L$, which implies that a performance quite similar to that of an approach without auxilliary variables can be achieved. \\
There exist many possible choices of constructing an initial guess for the $\mathbf{c}$-component of the iterative scheme. Common choices are algebraic methods, most notably transfinite interpolation \cite{gordon1973transfinite}. Once the $\mathbf{c}$-component has been computed with one of the available methods, a reasonable way to compute the corresponding $\mathbf{d}$-part is through a (separable) projection of $\mathbf{x}_{\xi}$ and $\mathbf{x}_{\eta}$ onto $\bar{\Sigma}$. \\
Slightly superior initial guesses can be generated using multigrid techniques as demonstrated in \cite{hinz2018elliptic}. The problem is first solved using a coarser basis and an algebraic initial guess, after which the coarse solution vector is prolonged and subsequently used as an initial guess. This is compatible with the techniques discussed in this section. However, instead of prolonging the full coarse solution vector, we only prolong the $\mathbf{c}$-component and compute the corresponding $\mathbf{d}$-component using an $L_2(\hat{\Omega})$-projection.

\subsection{Multipatch}
\label{subsect:Multipatch}
The reformulation with auxilliary variables has a particularly interesting application in multipatch-settings, especially when extraordinary patch vertices are present. Most of the techniques from subsection \ref{subsect:Solution_Single} are readily applicable but there exist subtle differences that shall be outlined in the following. \\
Let $\hat{\Omega}$ be a multipatch domain, i.e.,
\begin{align}
\hat{\Omega} = \bigcup \limits_{i = 1}^{n} \hat{\Omega}_i.
\end{align}
For convenience, let us assume that each $\hat{\Omega}_i$ is an affine transformation of the reference unit square $\tilde{\Omega} = [0, 1]^2$ with corresponding mapping $\mathbf{m}_i: \tilde{\Omega} \rightarrow \hat{\Omega}_i$, where
\begin{align}
\mathbf{m}_i(\mathbf{s}) = A_i \mathbf{s} + \mathbf{b}_i.
\end{align}
Here, $A_i$ is an invertible matrix, $\mathbf{b}_i \in \mathbb{R}^2$ some translation and the vector $\mathbf{s} = (s, t)^T$ contains the free variables in $\tilde{\Omega}$. The automated generation of a multipatch structure is a nontrivial task, which is not discussed in this paper. For an overview of possible segmentation techniques, we refer to \cite{buchegger2017planar, xiao2018computing, falini2019thb}. \\
Let $\tilde{\mathbf{x}}: \hat{\Omega} \rightarrow \Omega$ be such that $\tilde{\mathbf{x}}^{-1}: \Omega \rightarrow \hat{\Omega}$ is a harmonic mapping. Assuming that the $\hat{\Omega}_i$ are arranged such that $\hat{\Omega}$ is convex, Rado's theorem \cite{schoen1997lectures} applies and a harmonic $\tilde{\mathbf{x}}^{-1}$ is bijective. \\
In the case of a multipatch domain, pairs of faces $( \gamma^\alpha_i, \gamma^\beta_j ) \subset \partial \hat{\Omega}_i \times \partial \hat{\Omega}_j$ and sets of vertices $\{ \mathbf{p}_i, \ldots, \mathbf{p}_l \} \subset \partial \hat{\Omega}_i \times \ldots \times \partial \hat{\Omega}_l$  may coincide on $\hat{\Omega}$. As such, the bases $\Sigma$ and $\bar{\Sigma}$, whose elements constitute single-valued functions on $\hat{\Omega}$ are constructed from the patchwise discontinuous local bases $\Sigma_i$ and $\bar{\Sigma}_i$ with appropriate degree of freedom (DOF) coupling that canonically follows from the connectivity properties of the $\hat{\Omega}_i$.   In the multipatch case, we solve (\ref{eq:main_xi_eta}) by evaluating the associated integrals through a set of pull backs of the $\hat{\Omega}_i \subset \hat{\Omega}$ into the reference domain $\tilde{\Omega}$. Thanks to the affine nature of the pull back, replacement of $\boldsymbol{\xi}$-derivatives by local $\mathbf{s}$-derivatives is straightforward. \\
As such, the solution of (\ref{eq:main_xi_eta}) yields a collection of mappings $\{ \mathbf{x}_i \}_i$, with $\mathbf{x}_i: \tilde{\Omega} \rightarrow \Omega_i \subset \Omega$, where each $\mathbf{x}_i$ satisfies
\begin{align}
\label{eq:patch_function_composition}
\mathbf{x}_i \simeq \tilde{\mathbf{x}} \vert_{\hat{\Omega}_i} \circ \mathbf{m}_i.
\end{align}
As the right hand side of (\ref{eq:patch_function_composition}) is a composition of bijective mappings, the bijectivity of $\mathbf{x}_i$ depends on the quality of the approximation. If the $\mathbf{x}_i$ are bijective, they jointly form a parameterization of $\Omega$. \\
Unlike in the single-patch setting, the $L_2(\hat{\Omega})$-projection associated with the linear part of the residual vector is not separable. As such, the evaluation of vector products $A^{-1} B \mathbf{s}$ (see equation (\ref{eq:D_vec_approx})) becomes more involved. A possible workaround is explicit assembly and inversion of the Jacobian of the system (see equation (\ref{eq:Jacobian_full})), leading to increased computational times and memory requirements. \\
A possible alternative is the approximation of products of the form $A^{-1} B \mathbf{s}$ by a sequence of patchwise separable operations. In the following, we sketch a plausible approach. \\
Similar to the single-patch case, products of the form $(\mathbf{q}^u, \mathbf{q}^v)^T = A^{-1} B \mathbf{s}$ satisfy
\begin{align}
(\mathbf{q}^u, \mathbf{q}^v)^T & = \underset{ ( \tilde{\mathbf{q}}^u, \tilde{\mathbf{q}}^v )^T }{ \operatorname{argmin} } \sum_{i=1}^n \frac{1}{2} \int_{\hat{\Omega}_i} \left \| \begin{bmatrix} \mathbf{u}[ \tilde{\mathbf{q}}^u ] \\ \mathbf{v}[ \tilde{ \mathbf{q} }^v ] \end{bmatrix}  - \begin{bmatrix} \mathbf{x}^0_\xi[ \mathbf{s} ] \\  \mathbf{x}^0_\eta[ \mathbf{s} ] \end{bmatrix} \right \|^2 \mathrm{d} \boldsymbol{\xi}.
\end{align}
Let
\begin{align}
\tilde{\Sigma} = \bigcup \limits_{i=1}^{n} \bar{\Sigma}_i \equiv \{ \tilde{w}_i \}_i
\end{align}
be the patchwise discontinuous union of local (auxilliary variable) bases and let
\begin{align}
\tilde{ \mathbf{u} }[ \mathbf{g} ] & = \sum_i \mathbf{g}_i \tilde{w}_i, \nonumber \\
\tilde{ \mathbf{v} }[ \mathbf{h} ] & = \sum_i \mathbf{h}_i \tilde{w}_i.
\end{align}
In order to approximate $(\mathbf{q}^u, \mathbf{q}^v)^T$, we first find
\begin{align}
\label{eq:L2_multipatch_expanded}
(\mathbf{g}, \mathbf{h})^T = \underset{ ( \tilde{\mathbf{g}}, \tilde{\mathbf{h}} )^T }{ \operatorname{argmin} } \sum_{i=1}^n \frac{1}{2} \int_{\hat{\Omega}_i} \left \| \begin{bmatrix} \tilde{\mathbf{u}}[ \tilde{\mathbf{g}} ]  \\ \tilde{\mathbf{v}}[ \tilde{\mathbf{h}} ] \end{bmatrix} - \begin{bmatrix} \mathbf{x}^0_\xi[ \mathbf{s} ] \\  \mathbf{x}^0_\eta[ \mathbf{s} ] \end{bmatrix} \right \|^2 \mathrm{d} \boldsymbol{\xi}.
\end{align}
We perform a patchwise pullback of the $L_2$-projections into the reference domain where they are solved with the techniques from subsection \ref{subsect:Solution_Single}. Thanks to the affine nature of the pullback, the geometric factor associated with $\hat{\Omega}_i$ is constant and given by
\begin{align}
\det J_i = \det A_i.
\end{align}
Therefore, separability is not lost and the same efficiency as in the single-patch case is achieved. We restrict the solution of (\ref{eq:L2_multipatch_expanded}) to $\bar{\Sigma}$ by performing a weighted sum of components that coincide under coupling. Let $\bar{w}_i \in \bar{\Sigma}$ result from a coupling of $\{ \tilde{w}_\alpha, \ldots, \tilde{w}_\gamma \} \subset \tilde{\Sigma}$ and let $\{ \det J_{\alpha}, \ldots, \det J_ {\gamma} \}$ denote the set of corresponding local geometric factors. If the $\{ \tilde{w}_\alpha, \ldots, \tilde{w}_\gamma \} $ receive control points $\mathbf{g}_\alpha, \ldots, \mathbf{g}_\gamma$ under the projection, we set
\begin{align}
\label{eq:Multipatch_Restriction}
\mathbf{q}_i^u = \frac{ \det J_\alpha \mathbf{g}_\alpha + \ldots + \det J_\gamma \mathbf{g}_\gamma }{ \det J_\alpha + \ldots + \det J_\gamma },
\end{align}
and similarly for $\mathbf{q}^v$. Relation (\ref{eq:Multipatch_Restriction}) induces a canonical restriction operator from $\operatorname{span} \tilde{\Sigma}$ to $\operatorname{span} \bar{\Sigma}$ that is used to compute $(\mathbf{q}^u, \mathbf{q}^v)^T$ from $(\mathbf{g}, \mathbf{h})^T$.

\section{Numerical Experiments}
\label{sect:Numerical_Experiments}
In the following, we present several numerical experiments, demonstrating the functioning of the proposed algorithm. First, we present two single-patch problems after which we present a more involved multipatch parameterization. \\
In all cases, the auxiliary basis $\bar{\Sigma}$ results from one global $h$-refinement of the primal basis $\Sigma$.

\subsection{$L$-Bend}
\begin{figure}[h!]
\centering
\includegraphics[width=0.8 \textwidth]{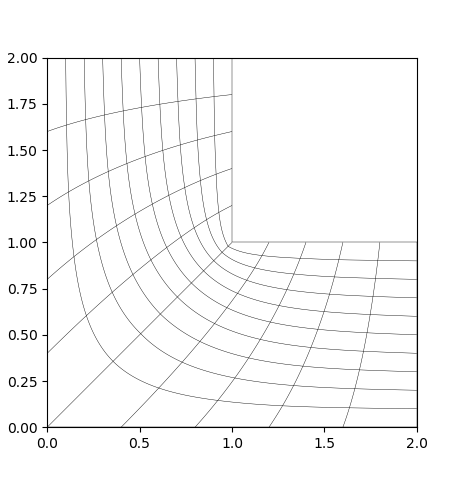}
\caption{Solution of the $L$-bend problem with the mixed-FEM algorithm.}
\label{fig:L_mixedFEM}
\end{figure}
\begin{figure}[h!]
\centering
\includegraphics[width=0.8 \textwidth]{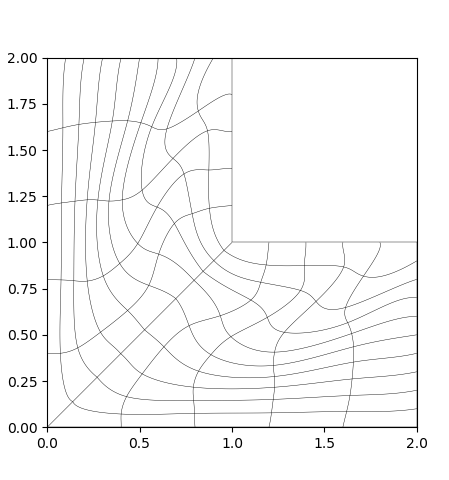}
\caption{Solution of the $L$-bend problem with constrained minimization of the \textit{Area Orthogonality} functional (see \cite{gravesen2012planar}).}
\label{fig:L_AO}
\end{figure}

\noindent As a proof of concept, we present results for the well-known single-patch $L$-bend problem. Wherever possible, we shall compare the results to a direct minimization of the Winslow-functional
\begin{align}
\label{eq:Winslow}
W( \mathbf{x} ) = \int_{ \hat{\Omega} } \frac{ g_{11} + g_{22} }{ \det J } \mathrm{d} \boldsymbol{\xi},
\end{align}
whose global minimizer (over $[\operatorname{span} \Sigma]^2$) coincides with a numerical approximation of the solution of (\ref{eq:EGG_not_discretized}) in the limit where $N \rightarrow \infty$ \cite{azarenok2009generation}. For the $L$-bend problem, we employ uniform cubic ($p=3$) knot-vectors in both directions with a $p$-fold knot-repetition at $\xi=0.5$ in order to properly resolve the $C^0$-continuity. As such we solve (\ref{eq:main_xi}) rather than (\ref{eq:main_xi_eta}). Figure \ref{fig:L_mixedFEM} shows the resulting parameterization along with the element boundaries under the mapping. The Schur-complement solver converges after $3$ iterations which amounts to $106$ evaluations of $R_N$. As can be seen in the figure, the parameterization is symmetric across the line connecting the upper and lower $C^0$-continuities which is expected behaviour from the shape of the geometry. We regard this as a positive sanity check for the functioning of the algorithm. Another observation is that despite the presence of knot-repetitions at $\xi=0.5$, the parameterization shows a large degree of smoothness along the corresponding isoline. Again, this is a positive result since the solution is expected to be an approximation of the global minimizer of (\ref{eq:Winslow}) (over $\mathbf{x} \in [\operatorname{span} \Sigma]^2$), which, in turn, approximates a smooth function. A substitution of the solution vector $\mathbf{c}_{\text{mf}}$ of the system of equations (\ref{eq:main_xi}) in (\ref{eq:Winslow}) gives
\begin{align}
W( \mathbf{c}_{\text{mf}} ) \simeq 3.01518,
\end{align}
whereas the global minimizer $\mathbf{c}_{\text{W}}$ of (\ref{eq:Winslow}) over the same basis yields
\begin{align}
W( \mathbf{c}_{\text{W}} ) \simeq 3.01425.
\end{align}
This constitutes another positive sanity check as the results are very close, while a substitution of the PDE-solution is slightly above the global minimum. As such, the PDE-solution comes with all the undesired characteristics of EGG-schemes such as the tendency to yield bundled / spread isolines near concave / convex corners. This does not occur in parameterizations based on the techniques of \cite{gravesen2012planar} (see figure \ref{fig:L_AO}). However, the $L$-bend example is rather contrived since a good parameterization is easily constructed with algebraic techniques. Here, the results only serve as a proof of concept.

\subsection{Tube-Like Shaped Geometry}
\label{subsect:Application_Separator}
\begin{figure}[h!]
\centering
\includegraphics[width=0.45 \textwidth]{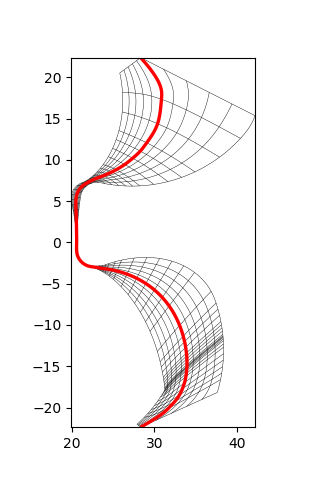} $\quad$
\includegraphics[width=0.45 \textwidth]{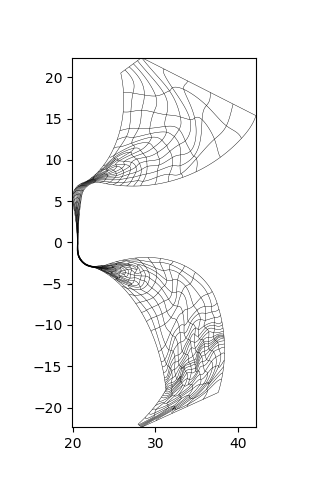}
\caption{PDE-based parameterization (left)  and area-orthogonality minimized parameterization (right) of a tube-like shaped geometry. }
\label{fig:separator_134}
\end{figure}

\begin{figure}[h!]
\centering
\includegraphics[width=0.45 \textwidth]{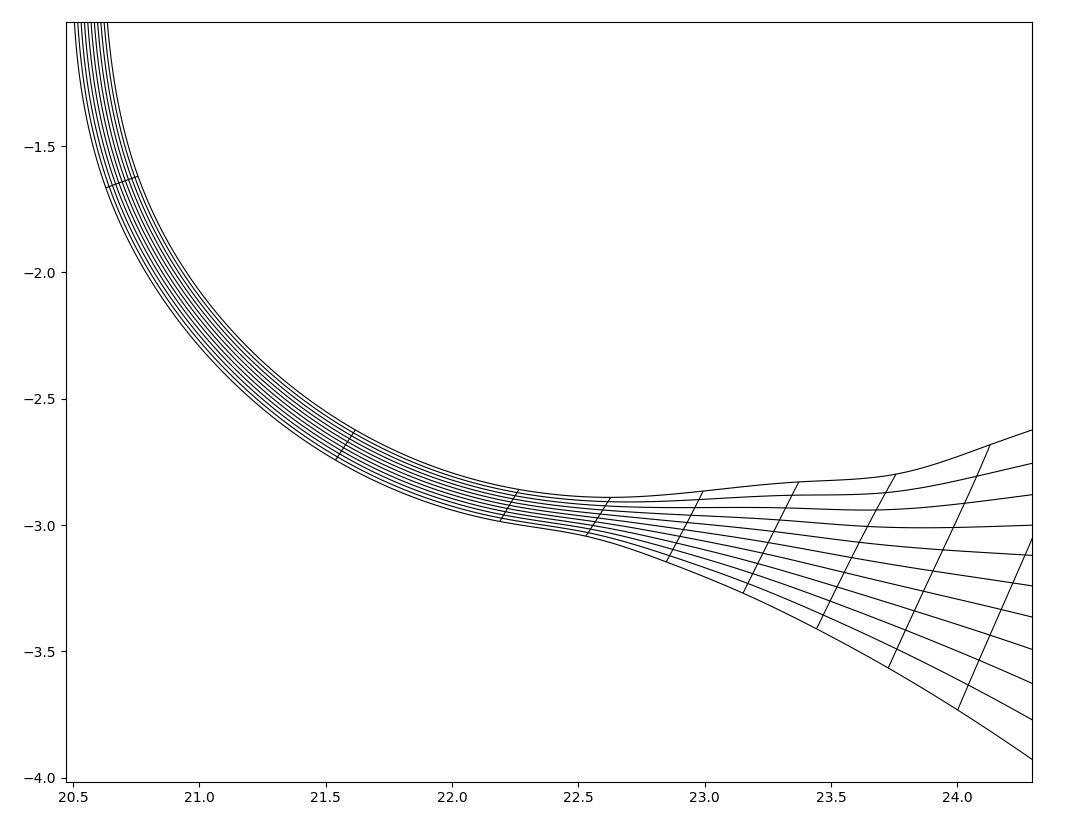} $\quad$
\includegraphics[width=0.45 \textwidth]{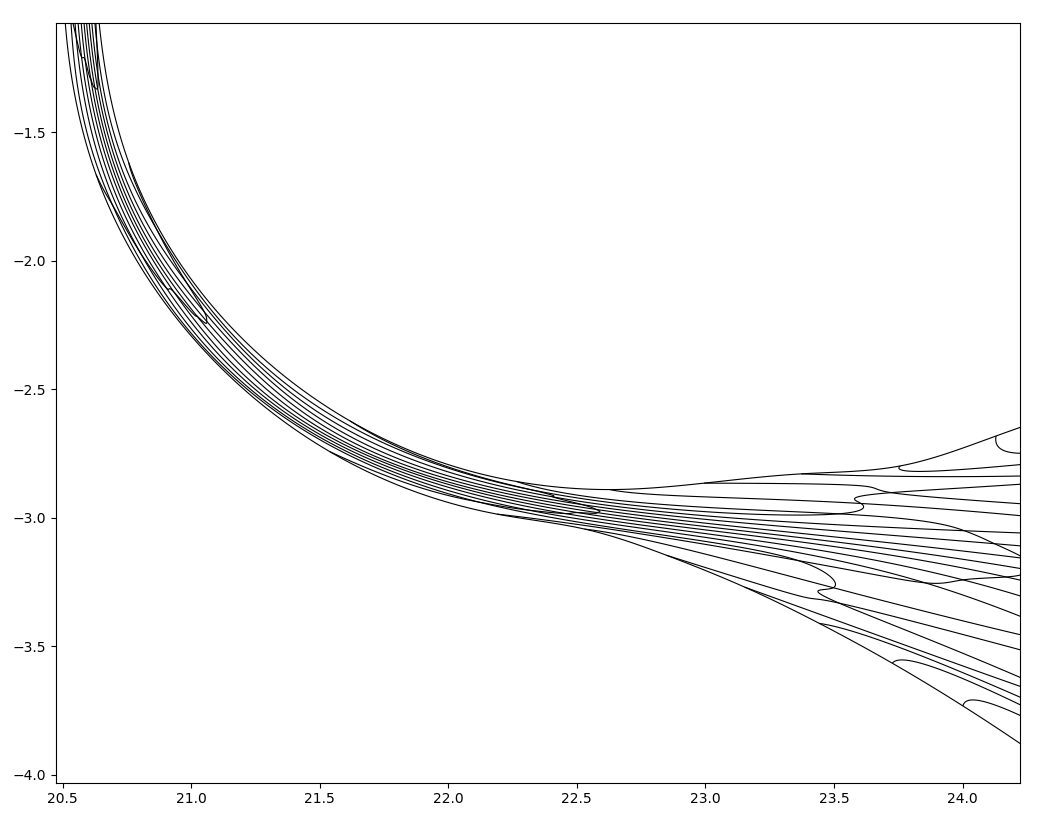}
\caption{Zoom-in on the PDE-based parameterization and area-orthogonality minimization parameterization.}
\label{fig:separator_134_zoom}
\end{figure}
\noindent In many cases, segmentation along knots with $p$-fold repetition and continuation with, for instance, techniques from \cite{hinz2018elliptic} on the smaller pieces is a viable choice. However, in some cases, a segmentation curve along which to split the geometry into smaller parts may be hard to find.  One such example is depicted in figure \ref{fig:separator_134} (left), which is a geometry taken from the practical application of numerically simulating a twin-screw machine. For convenience, the $\xi = 0.5$ isoline, across which the mapping is $C^0$-continuous, has been plotted in red. The usefullness of the proposed algorithm becomes apparent in this case: instead of having to generate a valid $\xi=0.5$ isoline, the isoline establishes itself from the solution of the PDE-problem. \\
As in the $L$-bend problem, we observe that the resulting parameterization exhibits a great degree of smoothness across the $\xi=0.5$ isoline, despite the continuity properties of $\Sigma$ and the spiked upper and lower boundaries. \\
The proposed algorithm produces superior results to the constrained optimization approach from \cite{gravesen2012planar} (see figure \ref{fig:separator_134}, right). In fact, here we initialized the optimization by the PDE-solution, as the solver struggles to find a feasible initial guess through optimization. This confirms the finding from \cite{hinz2018elliptic} that EGG-based approaches may be a viable alternative to finding feasible initial guesses for approaches based on optimization. Furthermore, we note the striking difference in the required number of iterations, which amount to over $100$ (constrained) in the optimization, while the PDE-solver converges in only $7$ iterations. \\
The poor performance of the optimization-approach can be explained by tiny gaps contained in the geometry, leading to natural jumps in the magnitude of the Jacobian determinant. As most cost functions are functions of the $g_{ij}$, they are very sensitive to jumps in $\det J$. This is further evidenced by the poor grid quality in the narrow part of the geometry (see figure \ref{fig:separator_134_zoom} right). In our experience, this is not the case for the PDE-solution (see figure \ref{fig:separator_134_zoom} left) and we successfully employed the approach for the automatic generation of a large number of similar geometries. \\
Finally, it should be noted that a comparison to the global minimizer of the Winslow-energy is not possible since the gradient-based optimizer we employed failed to further reduce the cost function from the evaluation of the PDE-solution.

\subsection{Multipatch Problem - The Bat Geometry}
\label{subsect:Bat}
\begin{figure}[h!]
\centering
\includegraphics[width=0.45 \textwidth]{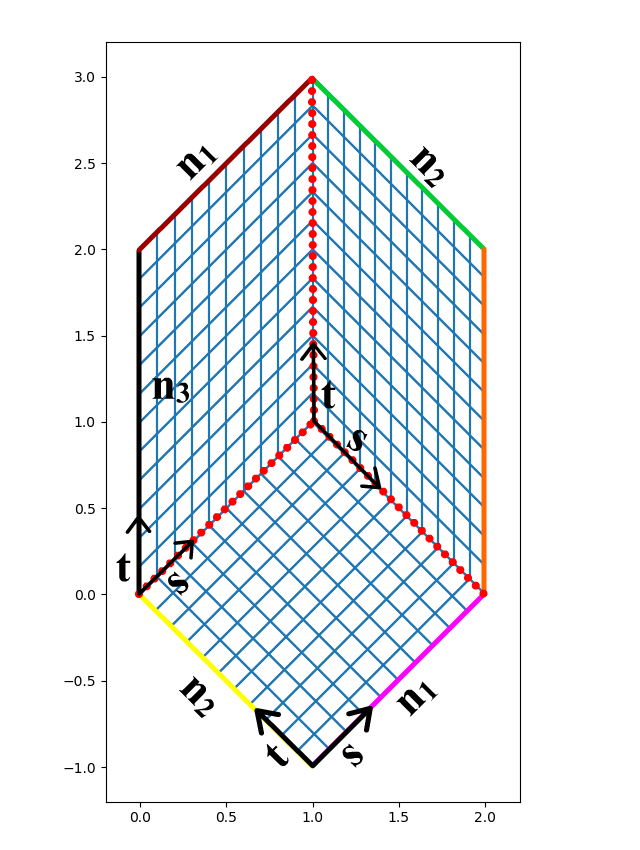} $\quad$ \includegraphics[width=0.45 \textwidth]{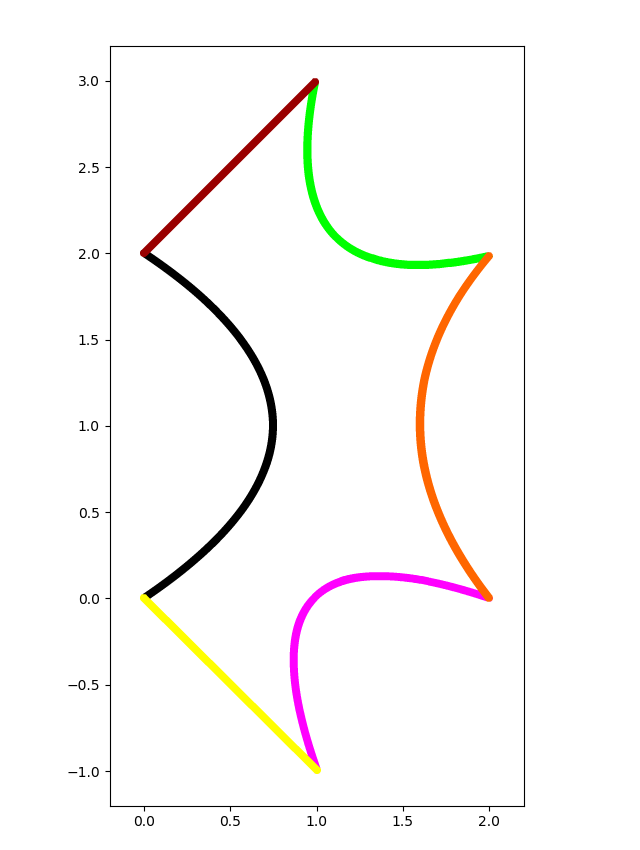}
\caption{Diamond shaped multipatch domain (left) and the target boundaries (right). Here, $n_1 = 10$, $n_2 = 11$ and $n_3 = 12$ denote the number of (uniformly-spaced) elements in each coordinate direction. There are no internal knot repetitions.}
\label{fig:diamond_domain}
\end{figure}

\begin{figure}[h!]
\centering
\includegraphics[width=0.45 \textwidth]{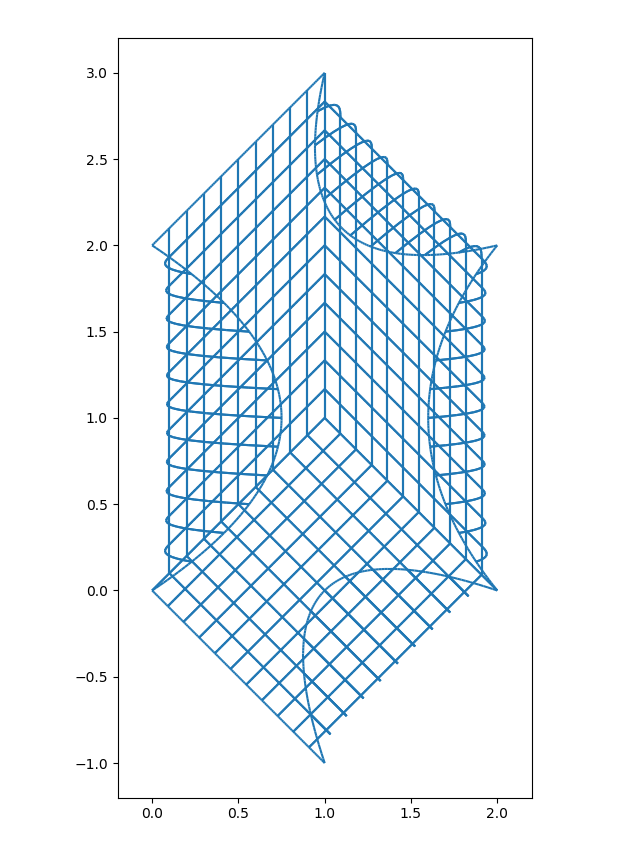} $\quad$ \includegraphics[width=0.45 \textwidth]{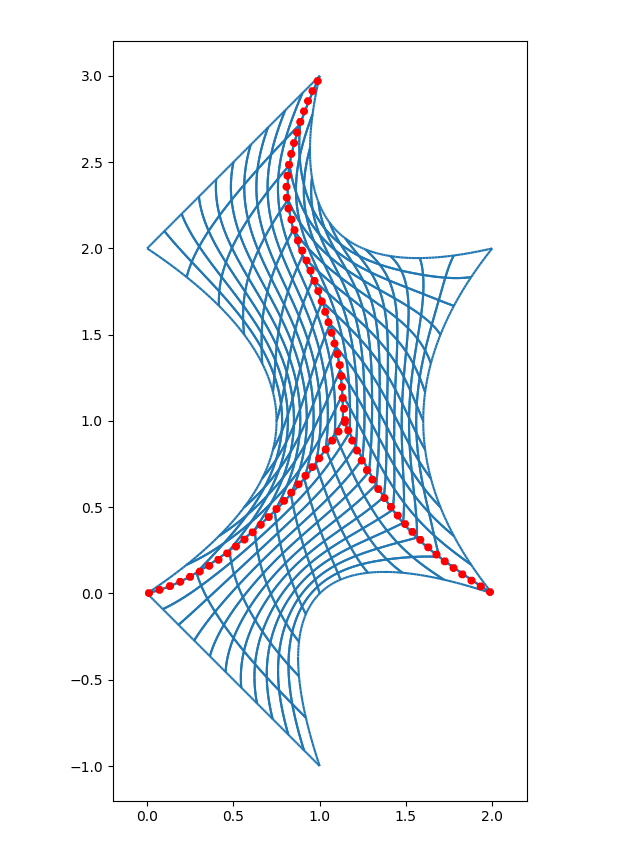}
\caption{The mapping that is passed on to the solver (left) and the resulting parameterization (right).}
\label{fig:init_result}
\end{figure}

Another interesting application of the proposed algorithm is that of a multipatch parameterization. In subsection \ref{subsect:Application_Separator}, we have successfully employed the algorithm to a geometry with a $C^0$-continuity along the $\xi=0.5$ isoline, which might as well be regarded as a two-patch parameterization with coupling along aforementioned isoline. A much more interesting multipatch application would be that of an uneven number of patches with extraordinary vertices. We are considering the diamond-shaped triple-patch domain depicted in figure \ref{fig:diamond_domain}, left. The target boundaries form the bat-shaped contour depicted in figure \ref{fig:diamond_domain}, right.  Note that, as required in subsection \ref{subsect:Multipatch}, the domain forms a convex subset of $\mathbb{R}^2$. For convenience we have highlighted the positions of the various boundaries under the mapping in different colors. Of course, of major interest shall be how the dotted red curve(s) in figure \ref{fig:diamond_domain} (left) are deformed under the mapping. 
Figure \ref{fig:init_result} (left) shows the mapping we utilize to initialize the Newton-Krylov solver while Figure \ref{fig:init_result} (right) shows the resulting geometry. Even though better initial guesses are easily constructed, here we have chosen to initialize the solver with a folded initial guess in oder to demonstrate that bijectivity is not a necessary condition for convergence. The Newton-Krylov solver converges after $6$ nonlinear iterations. The dotted red curves in figure \ref{fig:init_result} (right) show the internal interfaces of $\hat{\Omega}$ under the mapping.
\noindent We see that the patch interfaces are mapped into the interior of $\Omega$. The resulting geometry is bijective. However, the isolines make steep angles by the internal patch interfaces. This results from the additional pull back of $\tilde{\mathbf{x}} \vert_{\hat{\Omega}_i}$ into $\tilde{\Omega}$ via the operator $\mathbf{m}_i$ (see equation (\ref{eq:patch_function_composition})), which generally introduces a $C^0$-continuity in the composite mapping. Higher-order smoothness across patch interfaces is generally difficult to achieve and usually done by constructing bases whose elements possess higher-order continuity sufficiently far away from the extraordinary vertices. However, note that such bases may not allow for patchwise-affine transformations such that $L_2(\hat{\Omega}_i)$-projections lose their separability property. For a more rigorous definition of smoothness on multipatch topologies and strategies to build bases with local $C^{\geq1}$ smoothness on patch interfaces, we refer to \cite{buchegger2016adaptively}.

\section{Conclusion}
We have formulated an IgA-suitable EGG-algorithm that is compatible with spline bases $\Sigma$ possessing reduced regularity (whereby \textit{reduced} stands for global $C^{\geq 0}$-continuity) by introducing a set of auxilliary variables. We proposed an iterative Newton-Krylov approach operating on the Schur-complement of the linear part of the resulting nonlinear system of equations, which operates efficiently and reduces memory requirements. As such, it is suitable for large problems. Unlike similar $C^0$-compatible EGG-based approaches, the iterative solution method does not have to be initialized with a bijective mapping, significantly improving its usability in practice. However, this major advantage comes at the expense of increasing the problem size from $\simeq N$ to $\simeq N + c \vert \bar{\Sigma} \vert$, where $c=2$ or $c=4$, depending on the context. The impact is partially mitigated by the specially-taylored iterative solution algorithm. \\
We have presented three numerical experiments, two with a single patch and one resulting from a triple-patch configuration. In the single-patch case, we concluded that a substitution of the PDE-solution into the Winslow functional (equation (\ref{eq:Winslow})) yields an outcome that is close to that of the global function-minimizer (which is generally hard to find through direct minimization, due to the presence of $\det J$ in the denominator of equation (\ref{eq:Winslow})). As such, we concluded that the algorithm operates as expected and offers a viable alternative to direct minimization of (\ref{eq:Winslow}). However, it also comes with all the known drawbacks of EGG-based approaches and the two single-patch test cases demonstrate that it can yield inferior and superior results to other techniques, depending on the characteristics of the geometry. \\
As convergence is typically reached within only a few iterations, we conclude that the algorithm can serve as a computationally inexpensive method to initialize other methods that require a bijective initial guess. The required number of iterations can be further reduced by employing multigrid techniques (see \cite{hinz2018elliptic}) but this has not been implemented yet. \\
A major use case of the proposed algorithm is that of multipatch applications. In subsection \ref{subsect:Bat}, we presented results of the application to a triple-patch topology, where we successfully generated a patchwise bijective parameterization by approximating the composition of an inverse-harmonic mapping and patchwise affine transformations. The position of internal patch-interfaces under the mapping do not have to be imposed manually but follow naturally from the composite PDE-solution. \\
Finally, we observed that the composition with affine transformations results in nonsmooth transitions at patch interfaces. Higher-order smoothness can be achieved by a clever coupling of inter-patch DOFs sufficiently far away from extraordinary vertices.
\bibliographystyle{abbrv}
\bibliography{bibliography}
\end{document}